\documentclass[a4study,twoside,11pt]{article}
\usepackage{amssymb}
\usepackage{amsmath}
\usepackage{graphicx}
\usepackage{amsthm}
\usepackage{cite}
\usepackage{array}

\usepackage{amsmath,latexsym,amssymb,amsfonts,amsbsy,amsthm}
\usepackage{graphicx}
\usepackage{subfigure}
\usepackage[a4paper]{geometry}
\usepackage{epstopdf}
\usepackage{diagbox}
\usepackage{enumerate}
\usepackage{color}

\usepackage{comment}
\usepackage{indentfirst}

\newfont{\bb}{msbm10}

\newtheorem{example}{Example}[section]

\newtheorem{theorem}{Theorem}[section]

\newtheorem{lemma}{Lemma}[section]

\usepackage{verbatim}
\usepackage{algorithm}
\usepackage{algorithmicx}
\usepackage{algpseudocode}
\usepackage{amsmath}
\usepackage{graphics}
\usepackage{epsfig}
\usepackage{arydshln}

\usepackage{mathrsfs}
\usepackage{graphicx}
\usepackage{subfigure}
\usepackage{caption}
\usepackage{multirow}
\usepackage{threeparttable}
\numberwithin{equation}{section}

\usepackage{float}
\usepackage{footnote}
\makesavenoteenv{table}
\usepackage{booktabs}

\usepackage{authblk}

\usepackage{bm}
\usepackage{hyperref} 

\newenvironment{Proof}{{\noindent\it Proof.}}{\hfill $\square$ \par}

\title{ On multi-step extended maximum residual Kaczmarz method for solving large inconsistent linear systems} 

\author{
	A-Qin Xiao\\
	School of Mathematical Sciences, Tongji University,\\
	Shanghai, 200092, PR China.\\
	Email:xiaoaqin@tongji.edu.cn\\
	Jun-Feng Yin\thanks{Corresponding author.}\\
	School of Mathematical Sciences, Tongji University,\\
	Shanghai, 200092, PR China.\\
	Email:yinjf@tongji.edu.cn\\
	and\\
	Ning Zheng\\
	School of Mathematical Sciences, Tongji University,\\
	Shanghai, 200092, PR China.\\
	Email:nzheng@tongji.edu.cn\\
}

\date{ }
\begin{document}
 \cleardoublepage \pagestyle{myheadings}
\markboth{\small A.-Q. Xiao, J.-F. Yin and N. Zheng}
{\small On multi-step extended maximum residual Kaczmarz method }
\captionsetup[figure]{labelfont={bf},labelformat={default},labelsep=period,name={Fig.}}	
\captionsetup[table]{labelfont={bf},labelformat={default},labelsep=period,name={Table }}	

\maketitle

\begin{abstract}
A multi-step extended maximum residual Kaczmarz method is presented for the solution of the large inconsistent linear system of equations by using the multi-step iterations technique.
Theoretical analysis proves the proposed method is convergent and gives an upper bound on its convergence rate.
Numerical experiments show that the proposed method is effective and outperforms the existing extended Kaczmarz methods in terms of the number of iteration steps and the computational costs.\\



 \end{abstract}	

\section{Introduction}\label{sec1}
Consider the solution of the linear system of equations
\begin{equation}\label{eqincls}
Ax=b,
\end{equation}
where the matrix $A \in \mathbb{R}^{m \times n}$ and the vector $b \in \mathbb{R}^{m}$, 
which often arises from many practical scientific and engineering applications, for instance, image reconstruction \cite{2008HD}, signal processing \cite{2003B}, machine learning \cite{1991U} and option pricing \cite{2019FGNS}.
The Kaczmarz method is a simple and efficient iteration method for solving the linear system \cite{1937K}.
To improve the convergence of the Kaczmarz method, a randomized Kaczmarz method with an expected exponential convergence rate was proposed \cite{2009SV}.

However, the randomized Kaczmarz method fails to converge when the linear system \eqref{eqincls} is inconsistent.
To overcome the difficulty, by introducing an auxiliary vector $z_k$ to approximate $b_{\mathcal{R}(A)^{\perp}}$, i.e., the projection of $b$ onto the null space of $A^T$, the randomized extended Kaczmarz method was developed \cite{2013ZF}, which generates $z_{k+1}$ from the linear system $A^Tz=0$ and then computes the estimation solution $x_{k+1}$ from the linear system $Ax=b-z_k$. 
A tight upper bound for the convergence rate of the randomized extended Kaczmarz method was further presented \cite{2019D} by computing $x_{k+1}$ from $Ax=b-z_{k+1}$ since $z_{k+1}$ is a better approximation of $b_{\mathcal{R}(A)^{\perp}}$ than $z_{k}$. 
A cyclic column selection strategy was adopted for generating $z_{k+1}$ that makes the approximation of $b_{\mathcal{R}(A)^{\perp}}$ more precisely \cite{2019BW}.
For more studies on the extended Kaczmarz method, we refer the readers to \cite{2016PP,2020GLXQ,2022W,2023BW}.

In this paper, a multi-step extended maximum residual Kaczmarz method is presented by repeatedly implementing the iterate formula of $z_k$ many times to achieve a more accurate approximation of $b_{\mathcal{R}(A)^{\perp}}$. Moreover, the row index is determined by a greedy strategy to ensure the maximum entry of the residual vector is prioritized. 
Theoretical analysis gives an upper bound on the convergence rate of the proposed method.
Numerical experiments show that the proposed method is efficient and faster than existing methods.

The rest of this paper is organized as follows.
In Section \ref{secMEKs}, a multi-step extended greedy Kaczmarz method is proposed and its convergence theory is established. 
Numerical experiments are implemented to verify the efficiency of the proposed method in Section \ref{secnumer}.
Finally, some remarks and conclusions are drawn in Section \ref{secconclu}.

\section{The multi-step extended maximum residual Kaczmarz method} \label{secMEKs}
In this section, we first review the randomized extended Kaczmarz method and then propose a multi-step extended maximum residual Kaczmarz method for solving the inconsistent linear system of equations.

The randomized extended Kaczmarz method for the solution of the inconsistent linear system was first presented in \cite{2013ZF}, which iterates by two components
\begin{equation*} 
 z_{k+1}= z_{k}- \frac{ A_{(j)}^T z_{k} } {\lVert A_{(j)}\rVert^2_2}A_{(j)}\text{ and } x_{k+1}= x_{k}+\frac{b^{(i)}- A^{(i)}x_{k}-z^{(i)}_{k} }{\lVert A^{(i)}\rVert^2_2} (A^{(i)})^T,
\end{equation*}
where $A^{(i)}$ is the $i$th row of $A$, $A_{(j)}$ is the $j$th column of $A$, $b^{(i)}$ and $z^{(i)}$ are the $i$th entries of $b$ and auxiliary vector $z$, respectively.
Moreover, the expected convergence rates of iteration sequences $\{z_k\}_{k=0}^{\infty}$ and  $\{x_k\}_{k=0}^{\infty}$ generated by the randomized extended Kaczmarz method were derived, which are restated in Lemma \ref{lemconz} and Theorem \ref{thmREK}, respectively. 
\begin{lemma}{\rm\cite{2015NZZ} }\label{lemconz}    
\rm The sequence $\left\{ z_k\right\}_{k=0}^{\infty}$ generated in the iteration process of the randomized extended Kaczmarz method satisfies
  \begin{equation}\label{eqrekconvz}
      \mathbb{E}\left[ \left\|z_k-b_{\mathcal{R}(A)^{\perp} } \right\|_2^2\right]
      \leq \left(1-  \frac{\sigma_{\min}^2 (A)}{\left\| A\right\|_F^2} \right)^k\left\|  z_0-b_{\mathcal{R}(A)^{\perp}} \right\|_2^2,
  \end{equation}
where $b_{\mathcal{R}(A)^{\perp}}=b_{\mathcal{N}(A^T)}$ is the projection of $b$ onto the null space of $A^T$.
\end{lemma}

\begin{theorem}{\em\cite{2013ZF,2018BW}}\label{thmREK}
Let the initial vector $z_0=b$,  
the iteration sequence $\left\{ x_k\right\}_{k=0}^{\infty}$ generated by the randomized extended Kaczmarz method converges linearly in expectation to
the least squares solution $ x_{\ast}=A^{\dagger}b$.
Moreover, the solution error for the iteration sequence obeys 
\begin{equation*}
\mathbb{E}\left\| x_k- x_{\ast}\right\|_2^2 \leq  \alpha^{k-\lfloor\frac{k}{2} \rfloor}\left\| x_0- x_{\ast}\right\|_2^2+ \left(\alpha^{k-\lfloor\frac{k}{2} \rfloor} + \alpha^{ \lfloor\frac{k}{2} \rfloor} \right)\kappa^2(A)\|x_{\ast} \|_2^2, 
\end{equation*}
where $\alpha= 1- \frac{\sigma_{\min}^2 (A)}{\left\| A\right\|_F^2} $, $\lfloor k/2\rfloor$ is the floor of the constant $k/2$, 
$A^{\dagger}$, $\left\| A\right\|_F$, $\kappa(A)= {\sigma_{\max}(A)}/{\sigma_{\min}(A)}$, $\sigma_{\min}(A)$ and $\sigma_{\max}(A)$ are the Moore-Penrose inverse, Frobenius norm, condition number, smallest nonzero and largest singular value of $A$, respectively.
\end{theorem}
Note that the convergence rate of the sequence $\{x_k\}_{k=0}^{\infty}$ depends on the auxiliary sequence $\{z_k\}_{k=0}^{\infty}$.
To make $\{x_k\}_{k=0}^{\infty}$ converges to the least squares solution of system \eqref{eqincls} quickly,
we repeatedly execute the iterate formula of $z_k$ multiple times to obtain a more precise approximation of $b_{\mathcal{R}(A)^{\perp} }$ at each outer iteration. Moreover, the greedy strategy based on the maximum residual control is used to choosing the working row index.
For more studies on the greedy strategies, we refer the readers to \cite{2018BW,2020NZ,2021HM,2021ZL,2023XYZ}.   

In this work, a multi-step extended maximum residual Kaczmarz method is presented, which is described in Algorithm \ref{algMEKs}. 
When $\omega=1$, it gives the extended maximum residual Kaczmarz method.

\begin{algorithm}[!htbp]
 \caption{The multi-step extended maximum residual Kaczmarz method}\label{algMEKs}
   \begin{algorithmic}[1]
\Require $A,b,\ell,\omega,x_0\in \mathbb{R}^n$ and $z_0=b$
 \Ensure $x_\ell$
      \For {$k=0,1,\ldots,\ell-1$}
      \State Set $z_{0,k}=z_k$
     \For {$t=0,2,\ldots,\omega-1$} 
 \State Select $j_{t,k} \in \left\{ 1,2,\cdots, n \right\}$ with probability  $\text{Pr(column}= j_{t,k})= \frac{\|A_{(j_{t,k})}\|_2^2} { \|{A}\|_{{F}}^{2}} $
    \State Compute $z_{t+1,k}=z_{t,k}-\frac{ A_{(j_{t,k})}^T z_{t,k} }{\|A_{(j_{t,k})}\|_2^2}A_{(j_{t,k})}$
   \EndFor 
   \State Set $z_{k+1}= z_{\omega,k}$
        \State  Select $i_k=\arg\max\limits_{1\leq i\leq m}\Bigl| b^{(i)}- A^{(i)}x_{k}-z^{(i)}_{k+1}\Bigr|$ 
        \State Update  $ x_{k+1}= x_{k}+\frac{b^{(i_k)}- A^{(i_k)}x_{k}-z^{(i_k)}_{k+1} }{\lVert A^{(i_k)}\rVert^2_2} (A^{(i_k)})^T  $
 		\EndFor
 	\end{algorithmic}
 \end{algorithm}

Next, some lemmas are given before the convergence analysis of the multi-step extended maximum residual Kaczmarz method.
\begin{lemma}{\rm\cite{2018P}}\label{lemsqugeq} \rm
Let $\alpha_1, \beta_1$ be real numbers such that
  \begin{equation}\label{eqcond1}
      \alpha_1\in[0,1), \quad \beta_1\leq -1 \text{\quad and\quad} \beta_1-\alpha_1= \alpha_1\beta_1. 
  \end{equation} 
Then
\begin{equation}\label{eqcond2}
    (r_1+ r_2)^2\geq \alpha_1r_1^2-\beta_1r_2^2, \quad \forall r_1, r_2\in \mathbb{R}.
\end{equation}
\end{lemma}

\begin{lemma}{\rm\cite{2018P}}\label{lemalpbe}
\rm Let $\alpha_1, \beta_1$ satisfy \eqref{eqcond1}, then  
 \begin{equation}\label{eqcond3}
 \|x+y \|^2_2 \leq \alpha_1\| x\|^2_2- \beta_1\|y\|^2_2, \quad \forall x,y\in\mathbb{R}.
 \end{equation}
 \end{lemma}

\begin{lemma}\label{lemMEKconz}    
\rm The sequence $\{z_k \}_{k=0}^{\infty}$ generated in the iteration process of Algorithm \ref{algMEKs} converges to $b_{\mathcal{R}(A)^{\perp}}$, i.e., the projection of $b$ onto the null space of $A^T$ and satisfies
  \begin{equation}\label{eqconvz}
      \mathbb{E}\left[ \left\|z_{k+1}-b_{\mathcal{R}(A)^{\perp} } \right\|_2^2\right]
      \leq   \alpha ^{\omega (k+1)}\left\|  z_0-b_{\mathcal{R}(A)^{\perp} }  \right\|_2^2
       \leq \alpha ^{\omega (k+1)}\left\| b_{\mathcal{R}(A)} \right\|_2^2,
  \end{equation}
where $\alpha = 1- \frac{\sigma_{\min}^2 (A)}{\left\| A\right\|_F^2} $, the parameter $\omega\geq 1$ and $b_{\mathcal{R}(A)}$ is the range space of $A$.
 \end{lemma}

\begin{Proof}
    From the iteration formula of Algorithm \ref{algMEKs} 
    \begin{equation*}
        z_{t+1,k}= \left( I-\frac{A_{(j_{t,k})}A_{(j_{t,k})}^T}{\|A_{(j_{t,k})}\|_2^2} \right)z_{t,k}
    \end{equation*}
    and the fact
    \begin{equation*}
    \left( I-\frac{A_{(j_{t,k})}A_{(j_{t,k})}^T}{\|A_{(j_{t,k})}\|_2^2} \right)b_{\mathcal{R}(A)^{\perp}}=b_{\mathcal{R}(A)^{\perp}},
    \end{equation*}
it holds that
\begin{equation*}
     z_{t+1,k}-b_{\mathcal{R}(A)^{\perp}}= \left( I-\frac{A_{(j_{t,k})}A_{(j_{t,k})}^T}{\|A_{(j_{t,k})}\|_2^2} \right)( z_{t,k}-b_{\mathcal{R}(A)^{\perp}} ).
\end{equation*}
Then,
\begin{equation*}
     \begin{aligned}
    \lVert z_{t+1,k}-b_{\mathcal{R}(A)^{\perp}}\rVert^2_2= 
    \lVert z_{t,k}-b_{\mathcal{R}(A)^{\perp}}\rVert^2_2-  \frac{ |A_{(j_{t,k})}^T (z_{t,k}-b_{\mathcal{R}(A)^{\perp}})|^2 }{\|A_{(j_{t,k})}\|_2^2}.
    \end{aligned}
\end{equation*}
By taking the conditional expectation on the first $k$ outer iterations, it follows that 
 \begin{equation}\label{ineqconvz}
     \begin{aligned}
       \mathbb{E}_k\|z_{t+1,k}-b_{\mathcal{R}(A)^{\perp}} \|_2^2&= 
       \lVert z_{t,k}-b_{\mathcal{R}(A)^{\perp}}\rVert^2_2-  \sum_{j_{t,k}=1}^{n} \frac{\|A_{(j_{t,k})}\|_2^2} { \|{A}\|_{{F}}^{2} } \frac{ |A_{(j_{t,k})}^T (z_{t,k}-b_{\mathcal{R}(A)^{\perp}})|^2 }{\|A_{(j_{t,k})}\|_2^2}\\
       &=\lVert z_{t,k}-b_{\mathcal{R}(A)^{\perp}}\rVert^2_2- \frac{\|A^T (z_{t,k}-b_{\mathcal{R}(A)^{\perp}})\|_2^2}{\|{A}\|_{{F}}^{2} }    \\
      &\leq \left(1- \frac{\sigma_{\min}^2 (A)}{\left\| A\right\|_F^2}\right)\left\|  z_{t,k}-b_{\mathcal{R}(A)^{\perp} }  \right\|_2^2\\
      &\leq \alpha\left\|  z_{t,k}-b_{\mathcal{R}(A)^{\perp} } \right\|_2^2.
     \end{aligned}
 \end{equation}
Then, by using the full expectation on both side of the inequality \eqref{ineqconvz},  we have
\begin{equation}
\mathbb{E}\|z_{t+1,k}-b_{\mathcal{R}(A)^{\perp}} \|_2^2\leq \alpha \mathbb{E}\left\|z_{t,k}-b_{\mathcal{R}(A)^{\perp} } \right\|_2^2. 
\end{equation}
 Let $t=\omega-1$, since $z_{k+1}=z_{\omega,k}$, $z_{0,k}=z_k$, $z_0=b$ and $b= b_{\mathcal{R}(A)}+b_{\mathcal{R}(A)^{\perp}}$, it holds that
 \begin{equation} 
     \begin{aligned}
    \mathbb{E}\|z_{k+1}-b_{\mathcal{R}(A)^{\perp}} \|_2^2& \leq \alpha^{\omega} \mathbb{E}\|z_{k}-b_{\mathcal{R}(A)^{\perp}} \|_2^2\leq \alpha^{2\omega} \mathbb{E}\|z_{k-1}-b_{\mathcal{R}(A)^{\perp}} \|_2^2\leq\cdots \\
    &\leq \alpha^{\omega (k+1)} \|z_0-b_{\mathcal{R}(A)^{\perp}} \|_2^2\leq \alpha^{\omega (k+1)} \|b_{\mathcal{R}(A)} \|_2^2.
     \end{aligned}
 \end{equation}
Thus, we obtain the estimate \eqref{eqconvz}.
\end{Proof}
The convergence of the multi-step extended maximum residual Kaczmarz method is established in Theorem \ref{thmconvMEMRK}.

\begin{theorem}\label{thmconvMEMRK}
Let $\alpha_1,\beta_1$ be defined as in \eqref{eqcond1}, the sequence $\{ x^k\}_{k=0}^{\infty}$ generated by multi-step extended maximum residual Kaczmarz method converges to the least squares solution $ x_{\ast}=A^{\dagger}b$ of \eqref{eqincls}. Moreover, the solution error of the sequence satisfies
 \begin{equation}\label{eqconvMEK3}
 \left\|x_k-x_{\ast} \right\|_2^2 \leq \nu^k \left\|x_0-x_{\ast} \right\|_2^2+ \left( \nu^{k-\lfloor k/2\rfloor}+ \alpha^{\omega \lfloor k/2\rfloor} \right)\mu\gamma/{\alpha_1^2} \kappa^2(A) \|x_{\ast}\|_2^2,
  \end{equation}
where $\nu= 1- \frac{\alpha^2_1\sigma^2_{\min}(A)}{\gamma}$, $\alpha=1-\frac{\sigma^2_{\min}(A)}{\|A \|_F^2}$, $\mu= \frac{1+\beta_1}{ \min\limits_{1\leq i\leq m}\|A^{(i)}\|_2^2}+  \frac{\alpha_1\beta_1 }{\gamma}$ and $\gamma= m\max\limits_{1\leq i\leq m}\|A^{(i)}\|_2^2$.
\end{theorem}

\begin{Proof}
Substituting the decomposition $b= b_{\mathcal{R}(A)}+b_{\mathcal{R}(A)^{\perp}}$
 and $b= Ax_{\ast}$ into the iteration scheme of Algorithm \ref{algMEKs}, then
 \begin{equation}\label{eqinit}
     \begin{aligned}
x_{k+1}- x_{\ast}&=  x_k- x_{\ast}+ \frac{b_{\mathcal{R}(A)}^{(i_k)}-A^{(i_k)} x_k}{\|A^{(i_k)}\|_2^2}  (A^{(i_k)})^T + \frac{b_{\mathcal{R}(A)^{\perp}}^{(i_k)}-z_{k+1}^{(i_k)}  }{\|A^{(i_k)}\|_2^2}  (A^{(i_k)})^T \\
 &= \left(I_n- \frac{(A^{(i_k)})^TA^{(i_k)}}{\|A^{(i_k)}\|_2^2} \right)(x_k- x_{\ast}) + \frac{b_{\mathcal{R}(A)^{\perp}}^{(i_k)}-z_{k+1}^{(i_k)} }{\|A^{(i_k)}\|_2^2}  (A^{(i_k)})^T.
  \end{aligned}
\end{equation}
Since the two terms in the second equality of \eqref{eqinit} are perpendicular to each other, it follows that
 \begin{equation}\label{eqsumterms}
     \begin{aligned}
 \left\|x_{k+1}- x_{\ast}\right\|_2^2
 &= \left\|\left(I_n- \frac{(A^{(i_k)})^TA^{(i_k)}}{\|A^{(i_k)}\|_2^2} \right)(x_k- x_{\ast}) \right\|_2^2+ \left\|\frac{b_{\mathcal{R}(A)^{\perp}}^{(i_k)}-z_{k+1}^{(i_k)} }{\|A^{(i_k)}\|_2^2}  (A^{(i_k)})^T\right\|_2^2\\
     & =\|x_{k}-x_{\ast} \|_2^2- \frac{|A^{(i_k)}(x_k-x_{\ast}) |^2}{\lVert A^{(i_k)}\rVert^2_2}+ \frac{|b_{\mathcal{R}(A)^{\perp}}^{(i_k)}-z_{k+1}^{(i_k)} |^2}{\lVert A^{(i_k)}\rVert^2_2}\\
      &=\|x_{k}-x_{\ast} \|_2^2- \frac{|b_{\mathcal{R}(A)}^{(i_k)}-A^{(i_k)} x_k |^2}{\lVert A^{(i_k)}\rVert^2_2}+ \frac{|b_{\mathcal{R}(A)^{\perp}}^{(i_k)}-z_{k+1}^{(i_k)} |^2}{\lVert A^{(i_k)}\rVert^2_2}\\
      &=\|x_{k}-x_{\ast} \|_2^2- \frac{\Bigl|b_{\mathcal{R}(A)}^{(i_k)}+b_{\mathcal{R}(A)^{\perp}}^{(i_k)}-A^{(i_k)}x_{k}- z_{k+1}^{(i_k)}+ z_{k+1}^{(i_k)}-b_{\mathcal{R}(A)^{\perp}}^{(i_k)}\Bigr|^2 }{\|A^{(i_k)}\|_2^2}\\
      &\quad +\frac{|b_{\mathcal{R}(A)^{\perp}}^{(i_k)}-z_{k+1}^{(i_k)} |^2}{\lVert A^{(i_k)}\rVert^2_2} \\
      &\leq \|x_{k}-x_{\ast} \|_2^2- \frac{\alpha_1|b^{(i_k)}-A^{(i_k)} x_k- z_{k+1}^{(i_k)} |^2}{\lVert A^{(i_k)}\rVert^2_2}+ (1+\beta_1)\frac{|b_{\mathcal{R}(A)^{\perp}}^{(i_k)}-z_{k+1}^{(i_k)} |^2}{\lVert A^{(i_k)}\rVert^2_2}.
    \end{aligned}
\end{equation}
Here, the last inequality depends on Lemma \ref{lemsqugeq}.

In addition, from the choice $i_k=\arg\max\limits_{1\leq i\leq m}\Bigl|b^{(i)}-A^{(i)}x_k-z_{k+1}^{(i)} \Bigl|$ and the definition of residual $r_k= b-Ax_k-z_{k+1}$, it holds that 

\begin{equation}\label{eqmrstrategy}
\begin{aligned}
       \Bigl|b^{(i)}-A^{(i)}x_k-z_{k+1}^{(i)} \Bigl|^2
        &= \frac{\max\limits_{1\leq i\leq m}\Bigl|b^{(i)}-A^{(i)}x_k-z_{k+1}^{(i)} \Bigl|^2}{\lVert A^{(i_k)} \rVert^2_2\lVert r_k\rVert^2_2} \lVert A^{(i_k)}\rVert^2_2 \lVert r_k\rVert^2_2  \\
        & = \frac{\max\limits_{1\leq i\leq m}\Bigl|b^{(i)}-A^{(i)}x_k-z_{k+1}^{(i)} \Bigl|^2}{\sum\limits_{i=1}^m \Bigl|b^{(i)}-A^{(i)}x_k-z_{k+1}^{(i)} \Bigl|^2\lVert A^{(i_k)} \rVert^2_2}\lVert A^{(i_k)}\rVert^2_2 \lVert r_k\rVert^2_2  \\
        & \geq \frac{\lVert A^{(i_k)}\rVert^2_2}{ m\max\limits_{1\leq i\leq m} \lVert A^{(i)}\rVert^2_2 } \lVert (b_{\mathcal{R}(A)}-A x_k )+ (b_{\mathcal{R}(A)^{\perp}}-z_{k+1} ) \rVert^2_2\\
        & \geq \frac{\lVert A^{(i_k)}\rVert^2_2}{ m\max\limits_{1\leq i\leq m} \lVert A^{(i)}\rVert^2_2 } \lVert A(x_{\ast}-x_k)+ (b_{\mathcal{R}(A)^{\perp}}-z_{k+1} ) \rVert^2_2\\
        &\geq \lVert A^{(i_k)}\rVert^2_2\frac{\alpha_1\|A(x_k-x_{\ast}) \|_2^2-\beta_1\|b_{\mathcal{R}(A)^{\perp}}-z_{k+1} \|_2^2}{\gamma},
\end{aligned}
\end{equation}
where $\gamma= m\max\limits_{1\leq i\leq m} \lVert A^{(i)}\rVert^2_2$ and the last inequality holds because Lemma \ref{lemalpbe}. 

Combining \eqref{eqsumterms} and \eqref{eqmrstrategy}, then
\begin{equation}\label{eqsumtermsmek3}
    \begin{aligned}
     & \|x_{k+1}-x_{\ast} \|_2^2
      \leq \|x_{k}-x_{\ast} \|_2^2- \frac{\alpha_1|b^{(i_k)}-A^{(i_k)} x_k- z_{k+1}^{(i_k)} |^2}{\lVert A^{(i_k)}\rVert^2_2}+ (1+\beta_1)\frac{|b_{\mathcal{R}(A)^{\perp}}^{(i_k)}-z_{k+1}^{(i_k)} |^2}{\lVert A^{(i_k)}\rVert^2_2}\\
      &\leq \|x_{k}-x_{\ast} \|_2^2- \frac{\alpha_1^2\|A(x_k-x_{\ast}) \|^2_2}{\gamma}+ 
      \frac{\alpha_1\beta_1\|b_{\mathcal{R}(A)^{\perp}}-z_{k+1}\|^2_2}{\gamma}
      + \frac{(1+\beta_1)|b_{\mathcal{R}(A)^{\perp}}^{(i_k)}-z_{k+1}^{(i_k)} |^2}{\lVert A^{(i_k)}\rVert^2_2}\\
      &\leq \left(1- \frac{\alpha_1^2\sigma^2_{\min}(A)}{\gamma} \right)\|x_{k}-x_{\ast} \|_2^2
      +\left(\frac{\alpha_1\beta_1}{\gamma}+\frac{1+\beta_1}{\min\limits_{1\leq i\leq m}\|A^{(i)}\|_2^2} \right)\|b_{\mathcal{R}(A)^{\perp}}-z_{k+1} \|_2^2\\
      &\leq \nu\|x_{k}-x_{\ast} \|_2^2 +\mu\|b_{\mathcal{R}(A)^{\perp}}-z_{k+1} \|_2^2.
    \end{aligned}
\end{equation}
For integer $k\geq 0$, define $k_1= \lfloor k/2 \rfloor$ and $k_2=k-k_1$. Then, from \eqref{eqsumtermsmek3}, we obtain
\begin{equation}\label{eqestixk1}
    \begin{aligned}
  \left\|x_{k_1}-x_{\ast} \right\|_2^2 &\leq \nu  \left\|x_{k_1-1}-x_{\ast} \right\|_2^2+ \mu\|b_{\mathcal{R}(A)^{\perp}}- z_{k_1}\|_2^2\\
&\leq \nu^2 \left\|x_{k_1-2}-x_{\ast} \right\|_2^2+ \nu\mu \|b_{\mathcal{R}(A)^{\perp}}- z_{k_1-1}\|_2^2 + \mu \|b_{\mathcal{R}(A)^{\perp}}- z_{k_1}\|_2^2 \\
 &\leq \cdots\\
 &\leq \nu^{k_1}\left\|x_0-x_{\ast} \right\|_2^2+\sum_{l=0}^{{k_1}-1}\nu^l\mu\|b_{\mathcal{R}(A)^{\perp}}- z_{{k_1}-l}\|_2^2  \\
  &\leq   \nu^{k_1}\left\|x_0-x_{\ast} \right\|_2^2+\mu\|b_{\mathcal{R}(A)}\|_2^2 \sum_{l=0}^{{k_1}-1}\nu^l \\
   &\leq   \nu^{k_1}\left\|x_0-x_{\ast} \right\|_2^2+\mu\|b_{\mathcal{R}(A)}\|_2^2 \sum_{l=0}^{\infty}\nu^l \\
    &\leq   \nu^{k_1}\left\|x_0-x_{\ast} \right\|_2^2+ \frac{\mu\|b_{\mathcal{R}(A)}\|_2^2}{1-\nu}. 
    \end{aligned}
\end{equation}
The fifth inequality is obtained by Lemma \ref{lemMEKconz}.

From \eqref{eqsumtermsmek3} and Lemma \ref{lemMEKconz} again, it follows that
\begin{equation}\label{eqestix}
    \begin{aligned}
  \left\|x_k-x_{\ast} \right\|_2^2 & \leq \nu  \left\|x_{k_1+k_2-1}-x_{\ast} \right\|_2^2+ \mu  \|b_{\mathcal{R}(A)^{\perp}}- z_{k_1+k_2}\|_2^2   \\
 & \leq \cdots\\
 & \leq \nu^{k_2}\mathbb{E}\left\|x_{k_1}-x_{\ast} \right\|_2^2+\sum_{l=0}^{{k_2}-1}\nu^l\mu\|b_{\mathcal{R}(A)^{\perp}}- z_{{k_1+k_2}-l}\|_2^2\\
 & \leq \nu^{k_2}\mathbb{E}\left\|x_{k_1}-x_{\ast} \right\|_2^2+ \alpha^{\omega k_1}\sum_{l=0}^{{k_2}-1}\nu^l\mu\|b_{\mathcal{R}(A)^{\perp}}- z_{k_2-l}\|_2^2 \\
 & \leq \nu^{k_2}\mathbb{E}\left\|x_{k_1}-x_{\ast} \right\|_2^2+ \alpha^{\omega k_1} \mu\|b_{\mathcal{R}(A)}\|_2^2 \sum_{l=0}^{{k_2}-1}\nu^l\\
 & \leq \nu^{k_2}\mathbb{E}\left\|x_{k_1}-x_{\ast} \right\|_2^2 +\alpha^{\omega k_1}\frac{\mu\|b_{\mathcal{R}(A)}\|_2^2}{1-\nu}. 
    \end{aligned}
\end{equation}
Submitting \eqref{eqestixk1} into \eqref{eqestix}, then
\begin{equation}
    \begin{aligned}
   \left\|x_k-x_{\ast} \right\|_2^2 & \leq \nu^{k_2}\left( \nu^{k_1}\left\|x_0-x_{\ast} \right\|_2^2+ \frac{\mu\|b_{\mathcal{R}(A)}\|_2^2}{1-\nu} \right) + \alpha^{\omega k_1}\frac{\mu\|b_{\mathcal{R}(A)}\|_2^2}{1-\nu}\\
 & \leq \nu^k \left\|x_0-x_{\ast} \right\|_2^2+ \left( \nu^{k_2}+ \alpha^{\omega k_1} \right)\frac{\mu\gamma\sigma^2_{\max}(A)\|x_{\ast}\|_2^2}{\alpha_1^2\sigma^2_{\min}(A)}\\
  & \leq \nu^k \left\|x_0-x_{\ast} \right\|_2^2+ \left( \nu^{k_2}+ \alpha^{\omega k_1} \right)\mu\gamma\kappa^2(A)\|x_{\ast}\|_2^2 /{\alpha_1^2},
    \end{aligned}
\end{equation}
with $k_1=\lfloor k/2 \rfloor$ and $k_2= k-\lfloor k/2 \rfloor$. Here, the second inequality holds since $\|b_{\mathcal{R}(A)} \|_2^2=\|Ax_{\ast} \|_2^2\leq \sigma^2_{\max}(A)\|x_{\ast} \|_2^2$ and the definition of $\nu$.
\end{Proof}

\section{Numerical experiments} \label{secnumer}
In this section, numerical experiments are presented to verify the efficiency of the multi-step extended maximum residual Kaczmarz (MEMRK) method compared with the randomized extended Kaczmarz (REK) method \cite{2013ZF}, the partially randomized extended Kaczmarz (PREK) method  \cite{2019BW} and the extended maximum residual Kaczmarz (EMRK) method. 
Since the actual optimal parameter $\omega$ is difficult to obtain for the MEMRK method, we might choose the three parameters $\omega=4$ and $\omega=6$, named MEMRK1 and MEMRK2, respectively.

In the experiments, the right-hand side $b= Ax_{\ast}+\tilde{r}$ with the least squares solution $x_{\ast}=(1,1,\cdots,1)^T$ of the inconsistent linear system \eqref{eqincls}. Here, $\tilde{r}$ is a nonzero vector in the null space of $A^T$.
All iterations start from the initial vector $z_0=b$ and terminate when the relative residual vector (denoted as `RES') satisfies
$
{\rm RES}={\left\|r_k \right\|_2^2}/{\left\|r_0\right\|_2^2 } < 10^{-6}
$
or the number of iteration steps exceeds a maximum number, e.g., 50,000.

\begin{example}\em
The tested randomized matrices are generated with the standard normal distribution.
To guarantee the existence of vector $\tilde{r}$ for the underdetermined cases ($m\leq n$), the $m$-th row of $A$ is the average of its first two rows so that the matrix $A$ is rank-deficient.
\end{example}

In Tables \ref{tab:resultRandnov} and \ref{tab:resultRandnun},
the number of iteration steps (denoted as `IT') and the elapsed CPU time  in seconds (denoted as `CPU') for REK, PREK, EMRK and MEMRK methods with different $\omega$ are reported, respectively.

From Tables \ref{tab:resultRandnov} and \ref{tab:resultRandnun}, it is observed that all the extended Kaczmarz methods can successfully compute the solution of linear system with dense overdetermined or underdetermined coefficient matrix. 
Moreover, the MEMRK methods require fewer number of iteration steps and less elapsed CPU time than other extended Kaczmarz methods. It indicates that the multi-step strategy is efficient and can greatly improve the convergence of extended Kaczmarz methods.

\begin{figure}[!htbp]  
	\centering
	 \subfigure[$A\in \mathbb{R}^{6000\times 500}$]
	{
		\begin{minipage}[t]{0.48\linewidth}
			\centering
			\includegraphics[width=1\textwidth]{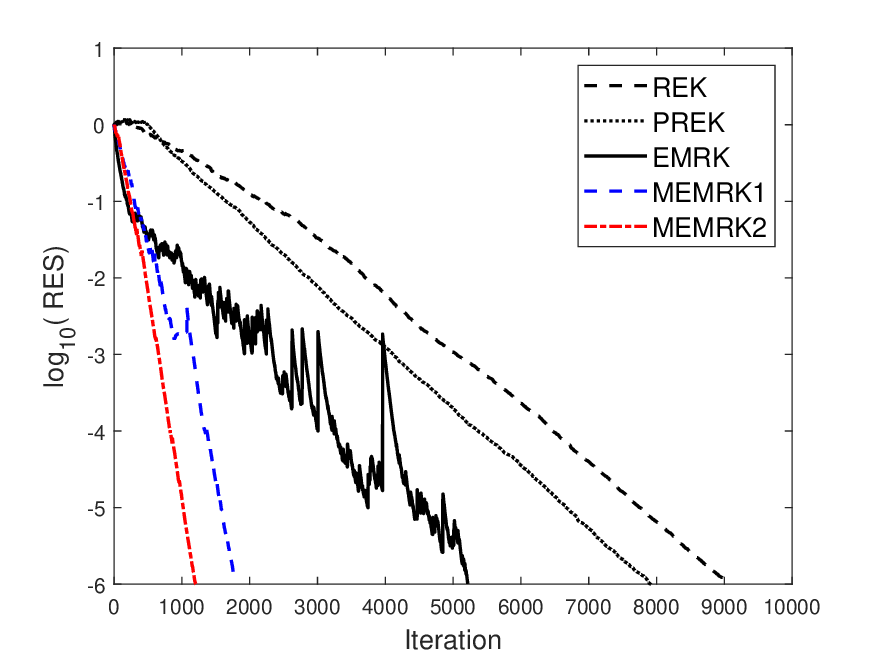}
		\end{minipage}
	} \hspace{-4mm} 
	 \subfigure[$A\in \mathbb{R}^{500\times 6000}$]
	{
		\begin{minipage}[t]{0.48\linewidth}
			\centering
			\includegraphics[width=1\textwidth]{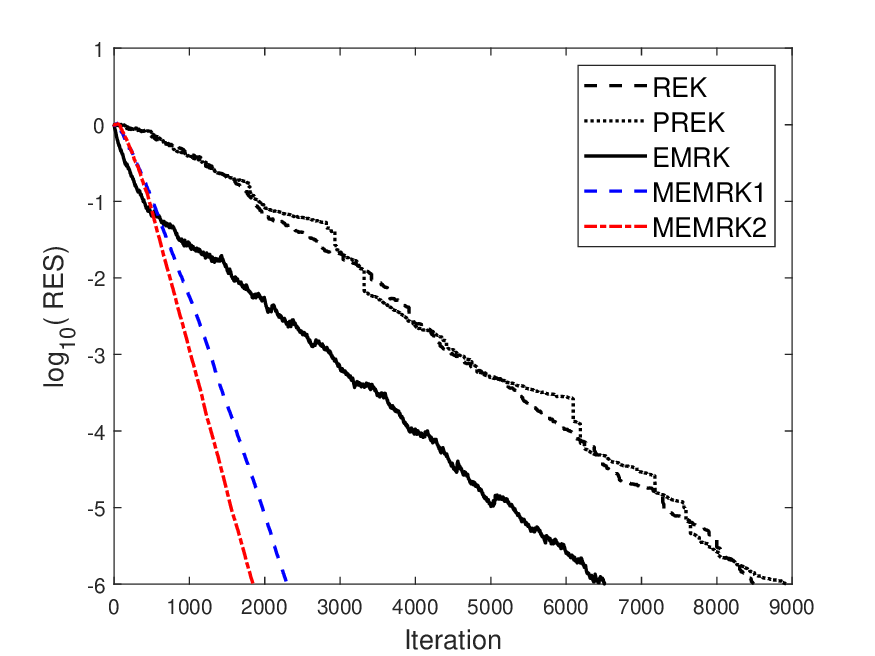}
		\end{minipage}
	}   
	\caption{ Convergence curves of the dense overdetermined (a) and underdetermined (b) cases.}
		\label{fig:RESvsITrandnovun}
\end{figure}

In Figure \ref{fig:RESvsITrandnovun}, the curves of the relative residual versus the number of iteration steps for all extended Kaczmarz methods are plotted, respectively. 

From Figure \ref{fig:RESvsITrandnovun}, it is observed that the curves of the MEMRK methods with different $\omega$ decrease much more faster than these of other extended Kaczmarz methods as the number of iteration steps increases in both the overdetermined and underdetermined cases, which shows the advantage of the multi-step strategy and verifies the numerical results in Tables \ref{tab:resultRandnov} and \ref{tab:resultRandnun}.

\begin{table}[!htbp] 
\scriptsize 
\centering 
\caption{Numerical results for overdetermined dense randomized matrices. }\label{tab:resultRandnov}  
\resizebox{\textwidth}{!}{  
\begin{tabular}{lllllll} 
\hline{Method}&{$ m \times n $} & 6000 $\times$ 500 & 7000 $\times$ 500 & 8000 $\times$ 500 & 9000 $\times$ 500 & 10000 $\times$ 500 \\ 
\hline  
  \multirow{2}{*}  {REK}              & IT&	 9084&	 9065&	 8899&	 8305&	 8460 \\  
& CPU&	 10.1820&	 10.7651&	 12.2949&	 13.5019&	 16.3095\\  
  \multirow{2}{*}  {PREK}             & IT&	 7913&	 8264&	 7721&	 7792&	 7707 \\  
& CPU&	 7.4323&	 9.3284&	 12.6688&	 11.9395&	 14.0162\\   
  \multirow{2}{*}  {EMRK} & IT&	 5216&	 5123&	 4674&	 4528&	 4657 \\  
& CPU&	 4.2567&	 5.0549&	 6.6650&	 6.1588&	 7.8184\\ 
  \multirow{2}{*}  {MEMRK1} & IT&	 1788&	 1622&	 1710&	 1584&	 1506 \\  
& CPU&	 1.9860&	 2.0718&	 2.5167&	 2.6749&	 2.8966\\  
  \multirow{2}{*}  {MEMRK2} & IT&	 1203&	 1343&	 1122&	 1151&	 1061 \\  
& CPU&	 1.5888&	 1.9689&	 1.9238&	 2.4030&	 2.2490\\  
\hline  
\end{tabular}  
 }   
\end{table} 

\begin{table}[!htbp] 
\scriptsize 
\centering 
\caption{Numerical results for  underdetermined dense randomized matrices. }\label{tab:resultRandnun}  
\resizebox{\textwidth}{!}{  
\begin{tabular}{lllllll} 
\hline{Method}&{$ m \times n $} & 500 $\times$ 6000 & 500 $\times$ 7000 & 500 $\times$ 8000 & 500 $\times$ 9000 & 500 $\times$ 10000 \\ 
\hline  
  \multirow{2}{*}  {REK}              & IT&	 8485&	 9062&	 8968&	 8706&	 8873 \\  
& CPU&	 8.8260&	 11.1547&	 13.3730&	 14.5079&	 17.3208\\  
  \multirow{2}{*}  {PREK}             & IT&	 8932&	 8513&	 8874&	 8233&	 7837 \\  
& CPU&	 8.5229&	 8.7046&	 10.7532&	 11.4435&	 12.5275\\   
  \multirow{2}{*}  {EMRK} & IT&	 6510&	 6430&	 6547&	 6168&	 6490 \\  
& CPU&	 6.7968&	 7.6881&	 9.2030&	 9.9112&	 12.3670\\ 
  \multirow{2}{*}  {MEMRK1} & IT&	 2294&	 2206&	 2263&	 2202&	 2191 \\  
& CPU&	 3.6826&	 4.1292&	 4.7821&	 5.1864&	 6.1458\\  
  \multirow{2}{*}  {MEMRK2} & IT&	 1844&	 1827&	 1736&	 1722&	 1686 \\  
& CPU&	 3.5415&	 4.0864&	 4.4443&	 4.8033&	 5.6269\\  
\hline  
\end{tabular}  
 }   
\end{table}

\newpage

\begin{example}\em
The tested sparse randomized matrices are obtained by the sparse standard normal distribution with a density of 0.1.
The $m$-th row of $A$ is the average of its first two rows for the underdetermined cases.
\end{example}

In Tables \ref{tab:resultspnov} and \ref{tab:resultspnun},
the number of iteration steps and the elapsed CPU time for REK, PREK, EMRK and MEMRK methods with different $\omega$ are listed, respectively.
 
From Tables \ref{tab:resultspnov} and \ref{tab:resultspnun}, it is seen that all the methods can converge to the solution of linear system with sparse overdetermined or underdetermined coefficient matrix. 
Moreover, the MEMRK methods need fewer number of iteration steps and less elapsed CPU time than other methods, which implies that the multi-step strategy is efficient and greatly improves the convergence. 

\begin{table}[!htbp] 
\scriptsize 
\centering 
\caption{Numerical results for overdetermined sparse randomized matrices.} \label{tab:resultspnov}  
\resizebox{\textwidth}{!}{  
\begin{tabular}{lllllll} 
\hline{Method}&{$ m \times n $} & 6000 $\times$ 1000 & 7000 $\times$ 1000 & 8000 $\times$ 1000 & 9000 $\times$ 1000 & 10000 $\times$ 1000 \\ 
\hline  
  \multirow{2}{*}  {REK}              & IT&	 22621&	 20670&	 20400&	 19528&	 19315 \\  
& CPU&	 22.3401&	 24.3057&	 27.9434&	 33.2451&	 36.6776\\  
  \multirow{2}{*}  {PREK}             & IT&	 18614&	 18098&	 17116&	 17039&	 16865 \\ 
& CPU&	 17.7720&	 20.9032&	 23.4030&	 27.7873&	 30.9092\\  
 \multirow{2}{*}  {EMRK} & IT&	 13974&	 12173&	 11217&	 12131&	 12145 \\  
& CPU&	 12.3488&	 12.8737&	 14.0063&	 19.4496&	 21.1520\\ 
  \multirow{2}{*}  {MEMRK1} & IT&	 4744&	 4717&	 4088&	 3953&	 3735 \\  
& CPU&	 4.9619&	 5.8414&	 5.8660&	 6.9518&	 8.4953\\  
  \multirow{2}{*}  {MEMRK2} & IT&	 3843&	 3250&	 3059&	 3043&	 3347 \\  
& CPU&	 4.5158&	 4.4045&	 4.7527&	 5.8405&	 8.1082\\  
\hline  
\end{tabular}  
 }   
\end{table}

\begin{table}[!htbp] 
\scriptsize 
\centering 
\caption{Numerical results for underdetermined sparse randomized matrices.} \label{tab:resultspnun}  
\resizebox{\textwidth}{!}{  
\begin{tabular}{lllllll} 
\hline{Method}&{$ m \times n $} & 1000 $\times$ 6000 & 1000 $\times$ 7000 & 1000 $\times$ 8000 & 1000 $\times$ 9000 & 1000 $\times$ 10000 \\ 
\hline  
  \multirow{2}{*}  {REK}              & IT&	 22034&	 21177&	 19645&	 19412&	 18890 \\  
& CPU&	 44.1519&	 47.6746&	 50.6162&	 56.9878&	 61.9422\\  
  \multirow{2}{*}  {PREK}             & IT&	 20421&	 20111&	 19196&	 18439&	 18972 \\  
& CPU&	 34.1488&	 40.1667&	 44.6941&	 48.9199&	 55.9808\\  
\multirow{2}{*}  {EMRK} & IT&	 14872&	 14091&	 13819&	 13649&	 13392 \\  
& CPU&	 28.7754&	 31.0587&	 35.0149&	 39.2104&	 43.0383\\
  \multirow{2}{*}  {MEMRK1} & IT&	 6044&	 5495&	 5210&	 4907&	 4670 \\  
& CPU&	 15.1553&	 15.5432&	 16.4901&	 17.3176&	 18.2764\\  
  \multirow{2}{*}  {MEMRK2} & IT&	 5070&	 4634&	 4202&	 4050&	 3930 \\  
& CPU&	 14.1914&	 14.3279&	 14.7004&	 15.7548&	 16.7525\\  
\hline  
\end{tabular}  
 }   
 \end{table} 

\begin{figure}[h]
	\centering
	 \subfigure[$A\in \mathbb{R}^{6000\times 500}$]
	{
		\begin{minipage}[t]{0.48\linewidth}
			\centering
			\includegraphics[width=1\textwidth]{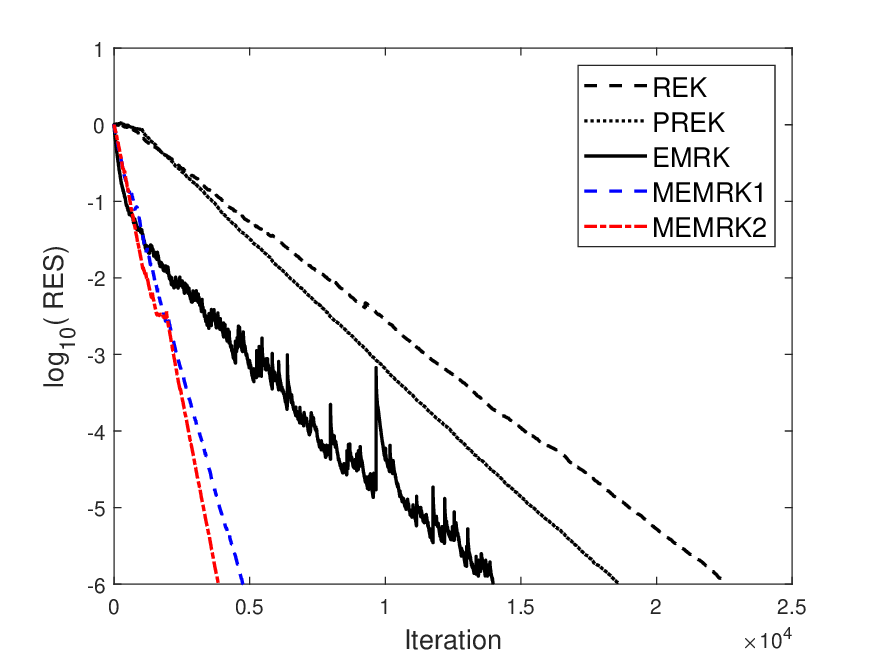}
		\end{minipage}
	} \hspace{-4mm} 
	 \subfigure[$A\in \mathbb{R}^{500\times 6000}$]
	{
		\begin{minipage}[t]{0.48\linewidth}
			\centering
			\includegraphics[width=1\textwidth]{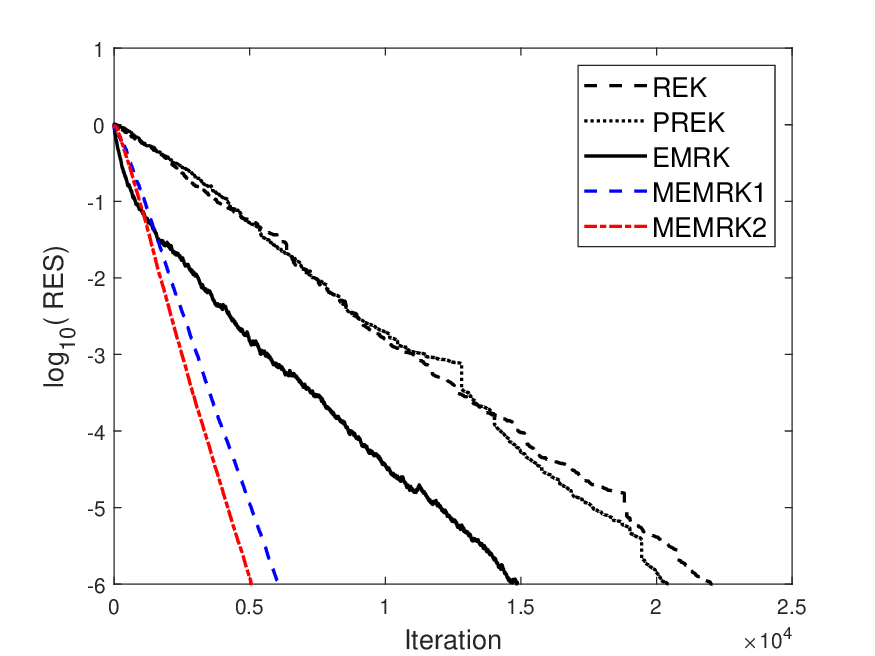}
		\end{minipage}
	}  
	\caption{ Convergence curves of the sparse overdetermined (a) and underdetermined (b) cases.}
		\label{fig:RESvsITspnovun}
\end{figure}

In Figure \ref{fig:RESvsITspnovun}, the curves of the relative residual versus the number of iteration steps for all methods are plotted, respectively.

From Figure \ref{fig:RESvsITspnovun}, it is observed that the convergence curves of MEMRK methods with different $\omega$ decrease much more faster than these of other methods as the number of iteration steps increases in both the overdetermined and underdetermined cases, which further confirms the numerical results in Tables \ref{tab:resultspnov} and \ref{tab:resultspnun}.

\begin{example} \em
The test example comes from the parallel-beam tomography medical image reconstruction problem generated by AIR Tools II \cite{2018HJ}.
In this experiment, the object domain is $[-20,20]\times [-20,20]$, angles are chosen with a step size of 2 from $0$ to $150$ and the distance between the first and last ray is 120. As a results, the size of the matrix $A$ is $9500 \times 1600$. The unique solution $x_{\ast}$ is obtained by reshaping a $40\times40$ original medical image.
The right-hand side vector is $b= Ax_{\ast} +\tilde{r}$, where $\tilde{r}$ is the Gaussian white noise vector with a noise level of $0.01$.
\end{example} 

The peak signal-to-noise ratio (PSNR)
$$
{\rm PSNR}(X_{t},X_{ r})=10 { \rm log}_{10} \frac{{ \max(X_{t}(i,j))}^2}{\frac{1}{mn}\sum\limits_{i=1}^{n}\sum\limits_{j=1}^{m}\left\|X_{t}(i,j)- X_{r}(i,j)\right\|^2 }
$$
is taken to measure the quality of the reconstruction results, where $X_{t}$ represents the true image of size $m\times n$, $X_{r}$ denotes the reconstructed image.
A larger PSNR value in dB indicates better preservation of the original image quality in the reconstructed image.

 \begin{figure}[!htbp] 
	\centering
        \subfigure[Exact phantom]
	{
		\begin{minipage}[t]{0.31\linewidth}
			\centering
			\includegraphics[width=1\textwidth]{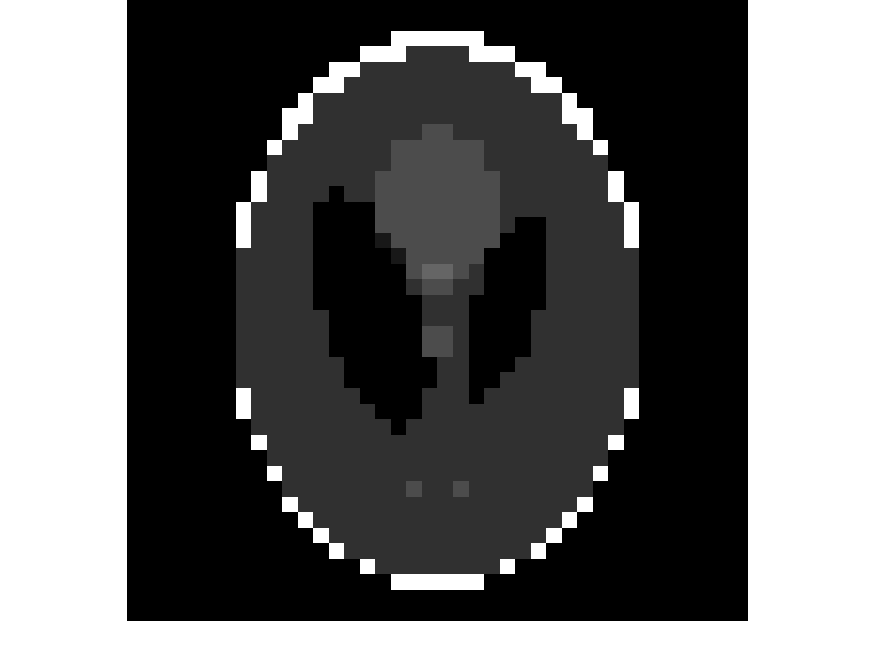} 
		\end{minipage}
	}\hspace{-14mm}
       \subfigure[\scriptsize REK]
	{
		\begin{minipage}[t]{0.33\linewidth}
			\centering
			\includegraphics[width=1\textwidth]{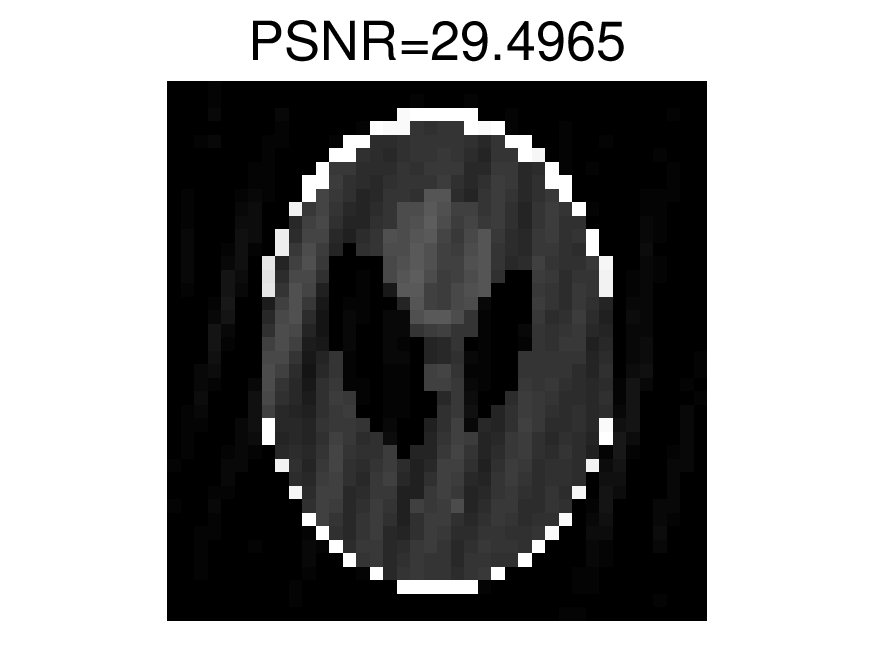} 
		\end{minipage}
	}\hspace{-14mm}
	\subfigure[\scriptsize PREK]
	{
		\begin{minipage}[t]{0.33\linewidth}
			\centering
			\includegraphics[width=1\textwidth]{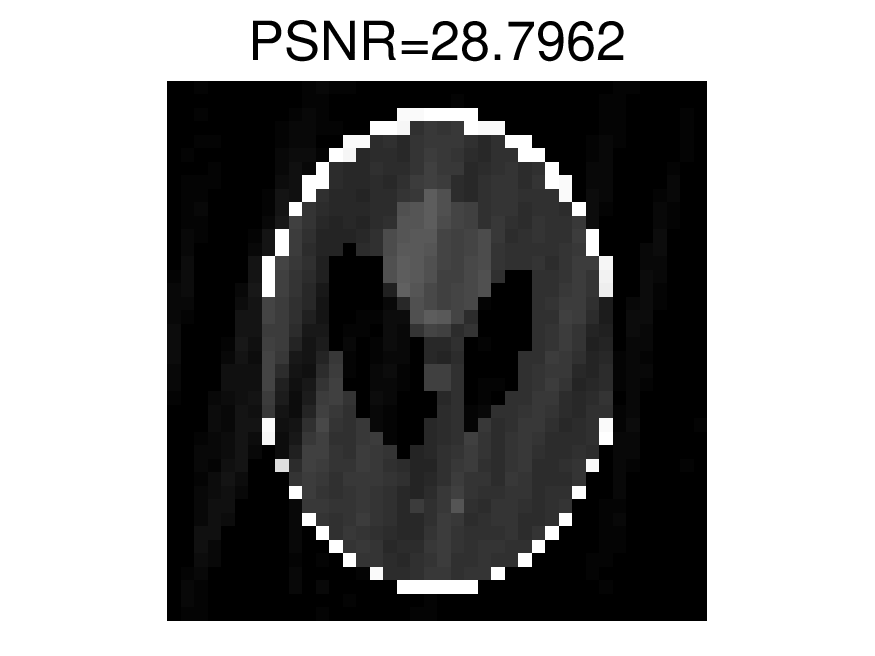} 
		\end{minipage}
	} \hspace{-14mm}
       \centering
	\subfigure[\scriptsize EMRK] 
	{
		\begin{minipage}[t]{0.33\linewidth}
			\centering
			\includegraphics[width=1\textwidth]{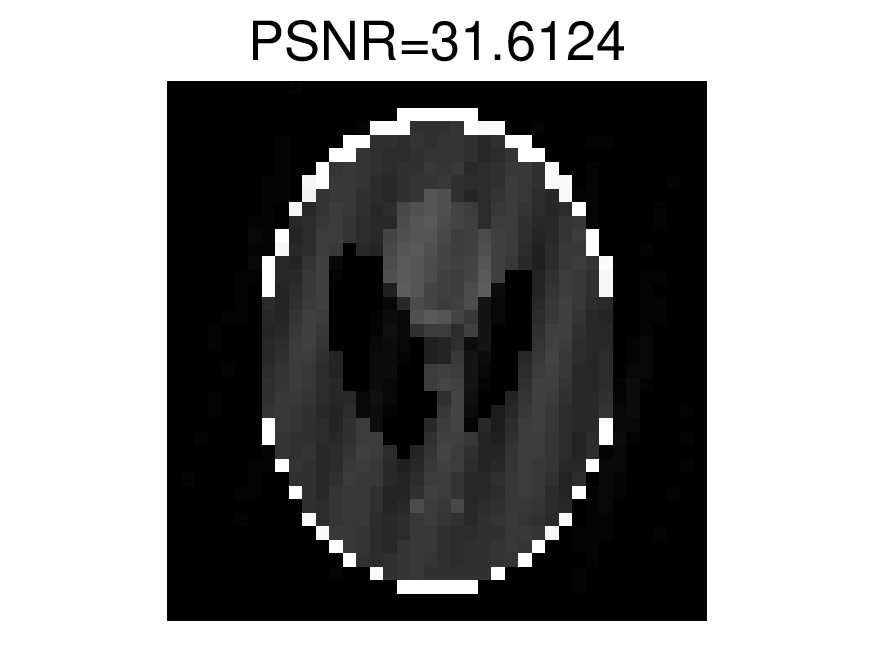} 
		\end{minipage}
	}\hspace{-14mm}
	\subfigure[\scriptsize MEMRK1] 
	{
		\begin{minipage}[t]{0.33\linewidth}
			\centering
			\includegraphics[width=1\textwidth]{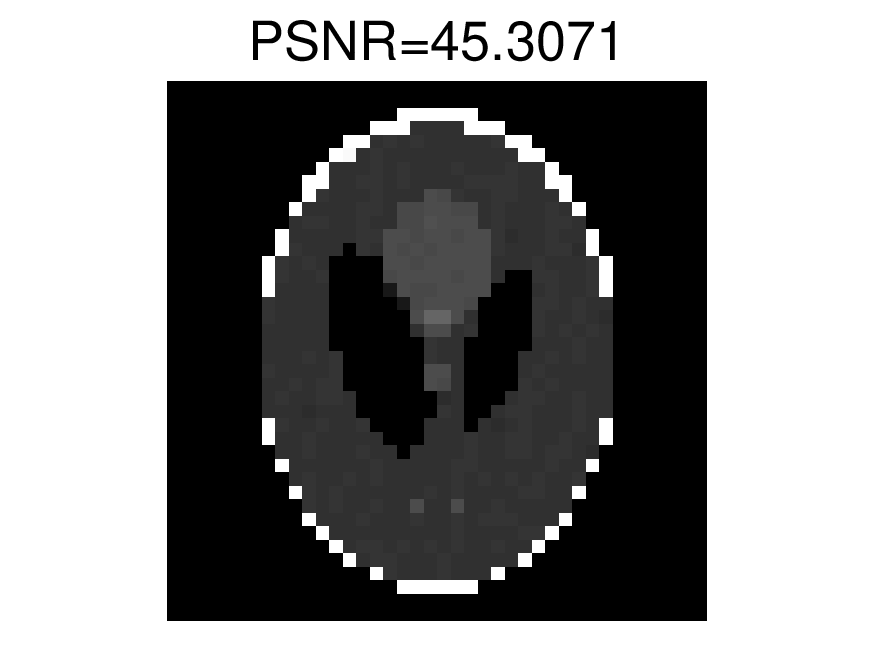} 
		\end{minipage}
	} \hspace{-14mm}
	\subfigure[\scriptsize MEMRK2]  
	{
		\begin{minipage}[t]{0.33\linewidth}
			\centering
			\includegraphics[width=1\textwidth]{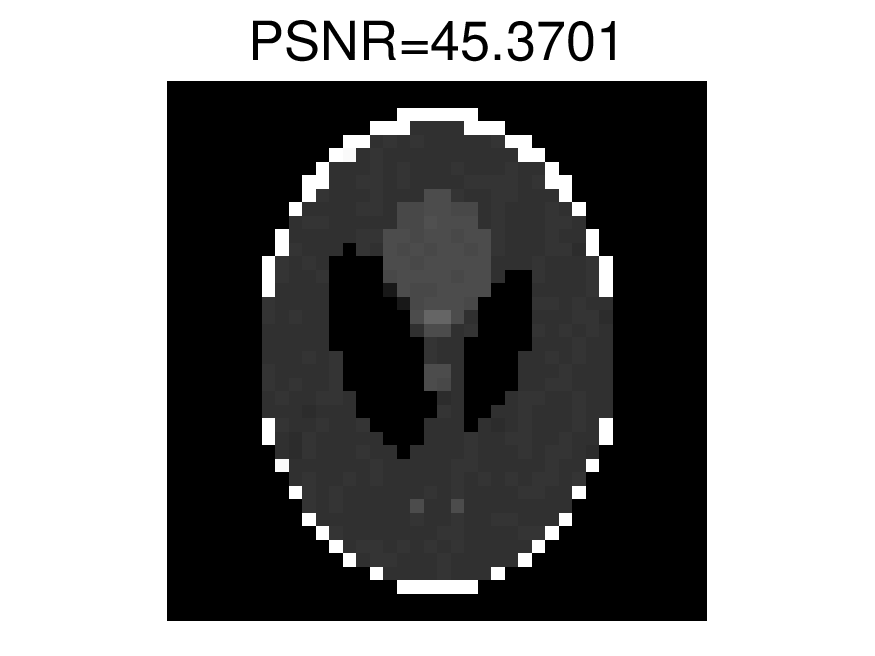} 
		\end{minipage}
	}
 \caption{Results of the parallel-beam tomography medical image problem.}
		\label{fig:pbtestpic}
\end{figure}

In Figure \ref{fig:pbtestpic}, the exact image and the recovered images obtained by  REK, PREK, EMRK and MEMRK methods with different $\omega$ are given after performing $10m$ iterations, where $m$ is the number of rows of $A$.

From Figure \ref{fig:pbtestpic}, it is observed that all methods can successfully recover the original image. The PSNR values of MEMRK methods are larger than those of other extended Kaczmarz methods. It implies that the images recovered by the MEMRK methods are much closer to the exact image than those obtained by other methods. 

\section{Conclusions}\label{secconclu}
In this paper, a multi-step extended maximum residual Kaczmarz method is developed for solving the inconsistent linear system of equations.
The convergence theory of the proposed method is established and the upper bound of the convergence rate for the method is derived.
Numerical experiments verify that the proposed method is efficient and superior to the existing extended Kaczmarz methods.  
 
\section*{Acknowledgements}
\noindent{\bf Funding} This work was supported by National Natural Science Foundation of China (No. 11971354). \\

\noindent{\bf Data availability statements } The datasets generated during the current study are available from the corresponding author on reasonable request.

\section*{Declarations}
\noindent{\bf Conflict of interest} The authors declare that they have no competing interests.

\bibliographystyle{plain}	 
\bibliography{refsmemrk.bib}

\begin{thebibliography}{10}

\bibitem{2023BW}
Zhong-Zhi Bai and Lu~Wang.
\newblock {On Multi-Step Randomized Extended Kaczmarz Method for Solving Large
  Sparse Inconsistent Linear Systems}.
\newblock {\em Applied Numerical Mathematics}, 2023.

\bibitem{2018BW}
Zhong-Zhi Bai and Wen-Ting Wu.
\newblock {On greedy randomized Kaczmarz method for solving large sparse linear
  systems}.
\newblock {\em SIAM Journal on Scientific Computing}, 40(1):A592--A606, 2018.

\bibitem{2019BW}
Zhong-Zhi Bai and Wen-Ting Wu.
\newblock {On partially randomized extended Kaczmarz method for solving large
  sparse overdetermined inconsistent linear systems}.
\newblock {\em Linear Algebra and Its Applications}, 578:225--250, 2019.

\bibitem{2003B}
Charles Byrne.
\newblock {A unified treatment of some iterative algorithms in signal
  processing and image reconstruction}.
\newblock {\em Inverse problems}, 20(1):103, 2003.

\bibitem{2019D}
Kui Du.
\newblock {Tight upper bounds for the convergence of the randomized extended
  Kaczmarz and Gauss--Seidel algorithms}.
\newblock {\em Numerical Linear Algebra with Applications}, 26(3):e2233, 2019.

\bibitem{2019FGNS}
Damir Filipovi{\'c}, Kathrin Glau, Yuji Nakatsukasa, and Francesco Statti.
\newblock {Weighted Monte Carlo with Least Squares and Randomized Extended
  Kaczmarz for Option Pricing}.
\newblock {\em Swiss Finance Institute Research Paper}, (19-54), 2019.

\bibitem{2020GLXQ}
Ying-Jun Guan, Wei-Guo Li, Li-Li Xing, and Tian-Tian Qiao.
\newblock {A note on convergence rate of randomized Kaczmarz method}.
\newblock {\em Calcolo}, 57:1--11, 2020.

\bibitem{2021HM}
Jamie Haddock and Anna Ma.
\newblock {Greed works: An improved analysis of sampling Kaczmarz--Motzkin}.
\newblock {\em SIAM Journal on Mathematics of Data Science}, 3(1):342--368,
  2021.

\bibitem{2018HJ}
Per~Christian Hansen and Jakob~Sauer J{\o}rgensen.
\newblock {AIR Tools II: algebraic iterative reconstruction methods, improved
  implementation}.
\newblock {\em Numerical Algorithms}, 79(1):107--137, 2018.

\bibitem{2008HD}
Gabor~T Herman and Ran Davidi.
\newblock {Image reconstruction from a small number of projections}.
\newblock {\em Inverse problems}, 24(4):045011, 2008.

\bibitem{1937K}
Stefan Karczmarz.
\newblock {Angenaherte auflosung von systemen linearer glei-chungen}.
\newblock {\em Bull. Int. Acad. Pol. Sic. Let., Cl. Sci. Math. Nat.}, pages
  355--357, 1937.

\bibitem{2015NZZ}
Deanna Needell, Ran Zhao, and Anastasios Zouzias.
\newblock {Randomized block Kaczmarz method with projection for solving least
  squares}.
\newblock {\em Linear Algebra and its Applications}, 484:322--343, 2015.

\bibitem{2020NZ}
Yu-Qi Niu and Bing Zheng.
\newblock {A greedy block Kaczmarz algorithm for solving large-scale linear
  systems}.
\newblock {\em Applied Mathematics Letters}, 104:106294, 2020.

\bibitem{2016PP}
Stefania Petra and Constantin Popa.
\newblock {Single projection Kaczmarz extended algorithms}.
\newblock {\em Numerical Algorithms}, 73:791--806, 2016.

\bibitem{2018P}
Constantin Popa.
\newblock {Convergence rates for Kaczmarz-type algorithms}.
\newblock {\em Numerical Algorithms}, 79(1):1--17, 2018.

\bibitem{2009SV}
Thomas Strohmer and Roman Vershynin.
\newblock {A randomized Kaczmarz algorithm with exponential convergence}.
\newblock {\em Journal of Fourier Analysis and Applications}, 15(2):262--278,
  2009.

\bibitem{1991U}
Shinji Umeyama.
\newblock {Least-squares estimation of transformation parameters between two
  point patterns}.
\newblock {\em IEEE Transactions on Pattern Analysis \& Machine Intelligence},
  13(04):376--380, 1991.

\bibitem{2022W}
Wen-Ting Wu.
\newblock {On two-subspace randomized extended Kaczmarz method for solving
  large linear least-squares problems}.
\newblock {\em Numerical Algorithms}, 89(1):1--31, 2022.

\bibitem{2023XYZ}
A-Qin Xiao, Jun-Feng Yin, and Ning Zheng.
\newblock {On fast greedy block Kaczmarz methods for solving large consistent
  linear systems}.
\newblock {\em Computational and Applied Mathematics}, 42(3):119, 2023.

\bibitem{2021ZL}
Yan-Jun Zhang and Han-Yu Li.
\newblock {Block sampling Kaczmarz--Motzkin methods for consistent linear
  systems}.
\newblock {\em Calcolo}, 58(3):39, 2021.

\bibitem{2013ZF}
Anastasios Zouzias and Nikolaos~M Freris.
\newblock {Randomized extended Kaczmarz for solving least squares}.
\newblock {\em SIAM Journal on Matrix Analysis and Applications},
  34(2):773--793, 2013.

\end{thebibliography}

\end{document}